\documentclass[a4paper,reqno,12pt]{amsart}
\usepackage{amsmath,amscd,amssymb,upgreek} 
\usepackage[all]{xy}
\usepackage{ifthen}
\usepackage{yhmath}
\usepackage{enumerate}
\usepackage[pdftex]{graphicx}
\usepackage[pdftex]{graphics}
\usepackage{hyperref}
\usepackage{pdfsync}
\usepackage{color}
\usepackage{hyperref}

\usepackage[applemac] {inputenc}

\newcommand{\application}[5]{
\ifthenelse{\equal{#1}{0}}{}{#1:}\begin{array}[t]{ccl}
	#2 & \longrightarrow & #3
	\ifthenelse{\equal{#4}{0}}{}{ \\ #4 & \longmapsto & #5} 
\end{array}}
\newcommand{\Exemple}[1]{\begin{Exple}\emph{#1}\end{Exple}}
\newcommand{\Remarque}[1]{\begin{Rem}\emph{#1}\end{Rem}}
\newcommand{\Theoreme}[1]{\begin{The}#1\end{The}}
\newcommand{\Proposition}[1]{\begin{Prop}#1\end{Prop}}
\newcommand{\Lemme}[1]{\begin{Lem}#1\end{Lem}}
\newcommand{\Corollaire}[1]{\begin{Cor}#1\end{Cor}}
\newcommand{\Definition}[1]{\begin{Def}#1\end{Def}}
\newcommand{\Preuve}[1]{\begin{proof}#1\end{proof}}

\def\ord{\operatorname{ord}}
\def\R{\mathbb{R}}
\def\C{\mathbb{C}}
\def\Pn{\mathbb{P}_{n}(\C)}
\def\H{\mathbb{H}}
\def\N{\mathbb{N}}
\def\Q{\mathbb{Q}}

\def\Z{\mathbb{Z}}

\def\S{\mathbb{S}}
\def\P{\mathbb{P}}

\def\Mc{\operatorname{\mathcal{M}c}}

\def\LVMB{\operatorname{LVMB}}
\def\LVM{\operatorname{LVM}}
\newcommand{\orb}{\operatorname{orb}}

\begin{document}
\newtheorem{The}{Theorem}[section]
\newtheorem{Prop}[The]{Proposition}
\newtheorem{Def}[The]{Definition}
\newtheorem{Lem}[The]{Lemma}
\newtheorem{Cor}[The]{Corollary}
\newtheorem{Rem}[The]{Remark}
\newtheorem{Exple}[The]{Example}

\title[A generalization of Sankaran and LVMB manifolds]{A generalization of Sankaran and LVMB manifolds \footnote{AMS classification 2010: 32J18, 32Q57, 32J99}}
\author[L. Battisti]{Laurent Battisti}
	\thanks{}
\author[K. Oeljeklaus]{Karl Oeljeklaus}
	\thanks{}
\begin{address}
	{
		\parbox{10cm}{
			\begin{tabular}{lcl}
				Laurent BATTISTI && Karl OELJEKLAUS \\ 
				RUHR-UNIVERSITÄT BOCHUM && AIX-MARSEILLE UNIVERSITÉ \\  
				FAKULTÄT FÜR MATHEMATIK && CNRS - LATP - UMR 7353\\
				Raum NA 4/26 && 39, Rue F. JOLIOT-CURIE\\
				D-44780 BOCHUM && F-13453 MARSEILLE\\
				Germany && France\\
				\href{mailto:laurent.battisti@rub.de}{\textbf{\texttt{laurent.battisti@rub.de}}}&&\href{mailto:karloelj@cmi.univ-mrs.fr}{\textbf{\texttt{karloelj@cmi.univ-mrs.fr}}}
			\end{tabular}
		}
	}
\end{address}

\begin{abstract}
In this paper we describe the construction of a new class of non-Kähler compact complex manifolds. They can be seen as a generalization of Sankaran, OT and LVMB manifolds. Moreover, we give properties of these new spaces. Their Kodaira dimension is $-\infty$ and under a mild condition they have algebraic dimension equal to zero.
\end{abstract}
\setcounter{section}{-1}
\maketitle
 
\section{Introduction}
\noindent In this article we construct a new family of non-Kähler complex compact manifolds by the combination of a method due to Bosio (\cite{Bosio:2001aa}) and another due to Sankaran (\cite{Sankaran:1989aa}). The class of manifolds constructed in this paper appears as a generalization of already known examples of non-Kähler manifolds, namely $\LVMB$ and Sankaran manifolds (ibid) along with OT manifolds (\cite{Oeljeklaus:2005aa}).\\

Although the field of non-Kähler geometry remains relatively unexplored, there is nevertheless constant progress. New classes of non-Kähler compact complex manifolds have been constructed and studied recently.

The first example is given by the class of $\LVMB$ manifolds. Their construction is due to Bosio (see \cite{Bosio:2001aa}) and can be summarized as follows. Given a family of subsets of $\{0, ..., n\}$ all having $2m+1$ elements (where $n$ and $m$ are integers such that $2m\leqslant n$) and a family of $n+1$ linear forms on $\C^{m}$ satisfying technical conditions, one can find an open subset $U$ of $\Pn$ and an action of a complex Lie group $G\cong\C^{m}$ on $\Pn$, such that the quotient $U/\C^{m}$ is a compact complex manifold. These manifolds generalize Hopf and Calabi-Eckmann manifolds and they also generalize a class of manifolds due to Meersseman (\cite{Meersseman:2000aa}), called $\LVM$ manifolds. \\

OT manifolds were constructed in \cite{Oeljeklaus:2005aa} by the second author and M.~Toma. One starts by choosing an algebraic number field $K$ having $s>0$ (resp. $2t>0$) real (resp. complex) embeddings. Then, for a nice choice of a subgroup $A$ of the groups of units $\mathcal{O}_{K}^{*}$ of $K$, the quotient $X(K,A)$ of $\H^{s}\times \C^{t}$ under the action of $A\ltimes \mathcal{O}_{K}$ is a complex compact manifold. The required condition on $A$ is that the projection on the $s$ first coordinates of its image through the logarithmic representation of units 
\begin{equation}\label{logUnit}\application{\ell}{\mathcal{O}_{K}^{*}}{\R^{s+t}}{a}{(\ln|\sigma_{1}(a)|, ..., \ln|\sigma_{s}(a)|, 2\ln|\sigma_{s+1}(a)|, ..., 2\ln|\sigma_{s+t}(a)|)}
\end{equation}
is a full lattice in $\R^{s}$. Such an $A$ is called \textbf{admissible}.
When there exists no proper intermediate field extension $\Q \subset K' \subset K$ with $A\subset \mathcal{O}_{K'}^{*}$, one says that the manifold $X(K,A)$ is of \textbf{simple type}.\\

In \cite{Sankaran:1989aa}, G.~Sankaran studies the action of a discrete group $W$ isomorphic to $\Z^{b}$ over the open subset $U$ of a toric manifold, chosen to obtain a compact quotient $U/W$. The group $W$ is a subgroup of the group of units of a number field $K$ and the construction of the infinite fan defining the toric manifold narrowly depends on $W$. This construction generalizes previous ones due to Inoue, Kato, Sankaran himself and Tsuchihashi. In the following, we will refer to the manifolds constructed in \cite{Sankaran:1989aa} as \emph{Sankaran manifolds}.\\

Our construction combines the methods of Bosio and Sankaran. First, we take the quotient of a well-chosen $n$-dimensional toric manifold $X_{\Delta}$ under the action of a complex Lie subgroup $G\cong\C^{t}$ of $(\C^{*})^{n}$ such that the quotient $X$ of $X_{\Delta}$ by $G$ is a (non-necessarily compact) manifold, then we choose a suitable open subset of $X$ on which a discrete group $W$ acts and gives a compact quotient. The group $W$ is a subgroup of the group of units of a number field $K$ having $s>0$ real and $2t>0$ complex embeddings, with $s+2t=n$. Of course, the choices of $W$ and $\Delta$ are strongly related, as in Sankaran's case. We obtain an $(s+t)$-dimensional non-Kähler complex compact manifold. \\

The paper is organized as follows: In section~\ref{PrelimGenSank} we introduce all notations and results that will be needed in the sequel. 

The next section is devoted to the description of the construction itself, which is subdivided into two steps. At the end of section~\ref{constr}, some remarks are made on the structure of the manifolds and we show how they generalize other known classes of manifolds.

In section~\ref{sousSectionInvariants} we provide results on invariants and geometric properties of the new manifolds. We will prove that their Kodaira dimension is $-\infty$, that they are not Kähler and finally show that their algebraic dimension is zero under a mild technical assumption. Moreover, we compute the second Betti number of OT manifolds under the same assumption, a fact which was already known for OT manifolds of simple type.

We conclude the paper by giving a concrete example of a $3$-dimensional manifold.


\section{\label{PrelimGenSank}Preliminaries}
\noindent Before starting the construction, we need to introduce or recall notations and results. Most of them can be found in Sankaran's paper (\cite{Sankaran:1989aa}), others come from \cite{Oeljeklaus:2005aa} and \cite{Battisti:2013aa}. Detailed proofs are only provided when modifications are needed to fit in with the current context.

The section is divided into four parts, as follows: First we settle the number-theoretic ground of our construction. Then, we study complex Lie groups which will arise later, in particular a certain Cousin group, i.e. a complex Lie group without non-constant holomorphic functions. Furthemore, we state some facts about manifolds with corners associated to a fan and finally construct a special fan.

\subsection{Number-theoretic notations and results}
Let $K$ be a number field having 
$n=s+2t$ distinct embeddings in $\C$ and let $\sigma_{1}, ..., \sigma_{s}$ be the $s$ real and $\sigma_{s+1}, ..., \sigma_{n}$ the $2t$ complex (non real) ones. Up to re-ordering these embeddings, we can assume that the relation $\sigma_{s+i} = \bar\sigma_{s+t+i}$ holds for $1 \leqslant i\leqslant t$. We also require that both $s$ and $t$ strictly positive.

Define the map \[\application{\sigma_{K}}{K}{\R^{s}\times \C^{2t}}{k}{(\sigma_{1}(k),..., \sigma_{s}(k),\sigma_{s+1}(k),..., \sigma_{n}(k)).}\]
Note $\mathcal{O}_{K}$ the set of algebraic integers of $K$. The image $\sigma_{K}(\mathcal{O}_{K})$ is a lattice of rank $n$ in $\R^{s}\times\C^{2t}\subset\C^{n}$. Let $\mathcal{O}_{K}^{*}$ be the group of units of $\mathcal{O}_{K}$. Since $s\geqslant1$, Dirichlet's Unit Theorem tells us that $\mathcal{O}_{K}^{*}$ is isomorphic to $\{\pm 1\}\times\Z^{s+t-1}$. The group $\mathcal{O}_{K}^{*}$ acts on $\C^{n}$ by componentwise multiplication: 
\[\eta\cdot(z_{1}, ..., z_{n}):= (\sigma_{1}(\eta) z_{1}, ..., \sigma_{n}(\eta) z_{n}).\]

For all $\eta \in \mathcal{O}_{K}^{*}$, set $\eta_{i}:=\sigma_{i}(\eta)$ if $i=1, ..., s$ and $\eta_{i}:=|\sigma_{i}(\eta)|$ if $i=s+1,..., s+t$. Call $\mathcal{O}^{*,+}_{K}$ the set of units $\eta \in \mathcal{O}_{K}^{*}$ such that $\eta_{i}>0$ for all $i=1,..., s$. The following theorem is due to Sankaran. For the convenience of the reader, we give the proof.

\Theoreme{[\cite{Sankaran:1989aa}, Theorem~2.1]\label{thSankExistenceW}For all $b \in \{1, ..., s\}$, there exists a subgroup $W < \mathcal{O}_{K}^{*,+}$ of rank $b$ satisfying the following condition, called ``\textbf{Assumption~C}''\footnote{In \cite{Sankaran:1989aa}, Sankaran formulates two others assumptions (A and B) which we do not use here.}: for all $\eta \neq 1 \in W$, either there is an $i \leqslant b$ such that $\eta_{i} > \eta_{j}$ for every $j>b$, or there is an $i\leqslant b$ such that $\eta_{i}<\eta_{j}$ for every $j>b$.}
\Preuve{Set $g:= s+ t-b$ and to every $\eta \in \mathcal{O}_{K}^{*,+}$ associate an element of $\R^{bg}$ via the map: 
\[\application{\varphi_{b}}{\mathcal{O}_{K}^{*,+}}{\R^{bg}}{\eta}{\left(\begin{array}{ccccc}
\ln(\eta_{1}/\eta_{b+1}) &  & \cdots & & \ln(\eta_{1}/\eta_{s+t})\\
 & & & & \\ 
\vdots& & \ln(\eta_{i}/\eta_{b+j}) & &\vdots \\
 & & & & \\ 
\ln(\eta_{b}/\eta_{b+1}) &  & \cdots & & \ln(\eta_{b}/\eta_{s+t})
\end{array}
\right).}\]
For $i = 1, ..., b$, let $\pi_{i,b}$ be the projection
\[\application{\pi_{i,b}}{\R^{bg}}{\R^{g}}{(t_{i,j})}{(t_{i,j})_{j=1,...,g}}\]
i.e. $\pi_{i,b}$ maps an element of $\R^{bg}$ written in matrix form to its $i$-th row.\\

The fact that $W$ satisfies Assumption~C means that for all $\eta \neq 1 \in W$, there exists $i \in \{1, ..., b\}$ such that \[\pi_{i,b}(\eta) \in Q_{b}:= \{(x_{1}, ..., x_{g})\in\R^{g}~|~ \text{the $x_{j}$'s are non-zero and of same sign}\}.\]

Call $U_{b}$ the vector space $\varphi_{b}(\mathcal{O}_{K}^{*,+}) \otimes \R \cong \R^{s+t-1} \subset \R^{bg}$. In order to prove the existence of a group $W$ of rank $b$ verifying Assumption~C, it is enough to find a $b$-dimensional linear subspace $A\subset U_{b}$ that is generated by elements of $\varphi_{b}(\mathcal{O}_{K}^{*,+})$ and such that for every $x\neq 0\in A$ there exists an $i \in\{1, ..., b\}$ with $\pi_{i,b}(x)\in Q_{b}$. We proceed by induction.\\

For $b=1$ the result is clear, as it is enough to take $A$ to be any line passing through the origin and any point of $Q_{1}$. 

Assume the result is proved at rank $b-1$: there exists a $(b-1)$-dimensional linear subspace $A'$ of $U_{b-1}$ having the following property: 
\[\text{for all }x \in A', \text{ there exists } i \in\{1, ..., b-1\} \text{ such that } \pi_{i,b-1}(x) \in Q_{b-1}.\]
In this case, for all $x \in A'$, the integer $i$ above satisfies: $\pi_{i,b}(x) \in Q_{b}$. \\
Now, fix a projection along $A'$, say $\pi_{A'}: U_{b} \rightarrow \R^{g}$. There is a map $L \in GL_{g}(\R)$ such that the following diagram commutes: 
\[\xymatrix{ & U_{b}\ar[ddl]_{\textstyle \smash{\pi_{b,b}}}\ar[ddr]^{\textstyle \smash{\pi_{A'}}} & \\
&&\\ \R^{g}\cong \pi_{b,b}(U_{b}) \ar@{-->}[rr]_{\textstyle L}& & U_{b}/A'\cong \R^{g}.}\]
Now we consider $Q_{A'}:= L(Q_{b})$. Take a line $\ell \subset Q_{A'}\cup \{0\}$, then the subspace $A:= \pi_{A'}^{-1}(\ell)$ fits our purpose. To see this, pick an $x \in A$. Either we have $x \in A'$ and in this case there is an $i \in \{1, ..., b\}$ such that $\pi_{i,b}(x) \in Q_{b}$, or $\pi_{A'}(x) \in \ell\setminus\{0\} \subset Q_{A'}$ and $\pi_{A'}(x) = L(\pi_{b,b}(x))$, which means that $\pi_{b,b}(x) \in Q_{b}$.}

\subsection{A Cousin group.}
Consider now the linear subspace \[H:=\{(0, ...,  0, z_{s+t+1}, ..., z_{s+t+t}) ~|~ z_{s+t+1}, ..., z_{s+t+t} \in \C\}\cong \C^{t}\] of $\C^{n}$ and call $\pi_{H}$ the projection from $\C^{n}$ to $\C^{s+t}$ with respect to $H$, given by the $s+t$ first coordinates.

\Lemme{\label{HcapKtrivial}The restriction of $\pi_{H}$ to $\sigma_{K}(\mathcal{O}_{K})$ is injective.}
\Preuve{It is sufficient to check that $H \cap \sigma_{K}(\mathcal{O}_{K}) = \{0\}$, which is straightforward because the $\sigma_{i}$'s are embeddings.}

Now let $W$ be as in the previous theorem for some $b\in\{1, ..., s\}$. Notice that the action of $W$ on $\C^{n}$ (resp. $\C^{n}/H\cong\C^{s+t}$) induces an action on $\sigma_{K}(\mathcal{O}_{K})$ (resp. $\pi_{H}(\sigma_{K}(\mathcal{O}_{K}))$).\\

Take an integral basis of $K$ and call $\mathcal{B}_{K}$ its image by $\sigma_{K}$. This is a basis of $\C^{n}$ over the complex numbers. Denote by $E$ the vector space generated over $\R$ by the vectors of $\sigma_{K}(\mathcal{O}_{K})$, note $\mathcal{B}:=(e_{1}, ..., e_{n})$ the canonical basis of $\C^{n}$ and $\mathcal{B}'$ the basis \[(e_{1}, ..., e_{s}, e_{s+1}+e_{s+t+1}, -i(e_{s+1}-e_{s+t+1}), ..., e_{s+t}+e_{s+2t}, -i(e_{s+t}-e_{s+2t})).\]
In the latter basis, the matrix of an element $\eta \in W$ is written as:
\[\left(\begin{array}{cccccr}
\sigma_{1}(\eta) &  &  & & & \raisebox{-10pt}{\text{\Huge{0~}}}\\
& \ddots & & & \\
& & \sigma_{s}(\eta) & & & \\
& & & \Re(\sigma_{s+t+1}(\eta)) & -\Im(\sigma_{s+t+1}(\eta)) & \\
& & & \Im(\sigma_{s+t+1}(\eta)) & \phantom{-}\Re(\sigma_{s+t+1}(\eta)) & \\
& & & & \ddots & \\
& & & & \multicolumn{1}{r}{\Re(\sigma_{n}(\eta))} & -\Im(\sigma_{n}(\eta)) \\
\smash{\text{\Huge{0}}} & & & & \multicolumn{1}{r}{\Im(\sigma_{n}(\eta))} & \phantom{-}\Re(\sigma_{n}(\eta)) 
\end{array}
\right).\] 

Moreover, $\mathcal{B}'$ is an $\R$-basis for the vector space $E$. To see this, it is enough to prove the following: 

\Lemme{\label{baseDeE}Let $P_{\mathcal{B}_{K},\mathcal{B}'}$ be the change-of-basis matrix from $\mathcal{B}_{K}$ to $\mathcal{B}'$. Then, all entries of $P_{\mathcal{B}_{K},\mathcal{B}'}$ are real, in other words one has $P_{\mathcal{B}_{K},\mathcal{B}'}\in GL_{n}(\R)$.}
\Preuve{Notice that $P_{\mathcal{B}_{K},\mathcal{B}'}=P_{\mathcal{B}_{K},\mathcal{B}}P_{\mathcal{B},\mathcal{B}'}$ and that the last $2t$ columns of $P_{\mathcal{B}_{K},\mathcal{B}}$ are pairwise conjugated, meanwhile the $s$ first columns are real. Hence, all the coefficients of the vectors $P_{\mathcal{B}_{K},\mathcal{B}} e_{j}$ for $j\in \{1, ..., s\}$ are real. Now, denote by $\overline{h_{1}}, ..., \overline{h_{t}}, h_{1}, ..., h_{t}$ the last $2t$ columns of $P_{\mathcal{B}_{K},\mathcal{B}}$ and observe that $P_{\mathcal{B}_{K},\mathcal{B}} (e_{s+j} + e_{s+t+j}) = 2\Re(h_{j})$ and $P_{\mathcal{B}_{K},\mathcal{B}} (-i(e_{s+j}-e_{s+t+j}))=-2\Im(h_{j})$ for $j=1,...,t$. The lemma is thus proven.}

Denote by $\widetilde{H}$ the vector space generated over $\R$ by the last $2t$ vectors of $\mathcal{B}'$. This is a $2t$-dimensional linear subspace of $E$ and we note $\pi_{\widetilde{H}} : E \rightarrow E/\widetilde{H}$ the quotient map.

\Remarque{\label{actionWDescend}The shape of the matrix of an element $\eta \in W$ in the basis $\mathcal{B}'$ shows that such an element induces a linear map from the real vector space $E\cong \R^{n}$ to itself and that this linear map descends to the quotient $E/ \widetilde{H} \cong \R^{s}$, where it is given by $\eta: (x_{1}, ..., x_{s}) \mapsto (\sigma_{1}(\eta)x_{1}, ..., \sigma_{s}(\eta)x_{s})$. Furthermore, one has $\widetilde{H} \cap \sigma_{K}(\mathcal{O}_{K}) = \{0\}$.}

Now we describe the linear subspace $H$ in the basis $\mathcal{B}_{K}$; it is generated by the $t$ vectors $(h_{1}, ..., h_{t})$, where $h_{i}$ is the vector $e_{s+t+i}$ of the canonical basis written in the basis $\mathcal{B}_{K}$ (according to notations of lemma~\ref{baseDeE}).~\\

The group $H$ is a closed Lie subgroup of $\C^{n}/\sigma_{K}(\mathcal{O}_{K})\cong (\C^{*})^{n}$, via the following map $\iota$: 

\begin{equation}\label{injectionCt}
\begin{array}{ccc}
H & \longrightarrow & (\C^{*})^{n} \\
& & \\
\left(\begin{array}{c}
0\\
\vdots\\
0\\
z_{s+t+1}\\
\vdots\\
z_{s+t+t}
\end{array}\right) & \longmapsto & 
\left(\begin{array}{c}
\exp(2i\pi \sum\limits_{i=1}^{t}h_{i,1}z_{s+t+i})\\
\\
\vdots\\
\\
\exp(2i\pi \sum\limits_{i=1}^{t}h_{i,n}z_{s+t+i})
\end{array}\right),
\end{array}
\end{equation}
where $h_{i,j}$ is the $j$-th coordinate of the vector $h_{i}$. We still denote by $H$ the image of $H$ under this map. We shall see that the group $(\C^{n}/\sigma_{K}(\mathcal{O}_{K}))/H$ only has trivial holomorphic functions.

A connected complex Lie group admitting no non-constant holomorphic functions is called a \emph{Cousin group}, or a \emph{toroidal group}. One has:

\Corollaire{\label{CstGpeCousin}The quotient of the complex Lie group $\C^{n}/\sigma_{K}(\mathcal{O}_{K})\cong(\C^{*})^{n}$ by $H$ is a Cousin group, which we call $C_{0}$. We denote by $p$ the quotient map $p: \C^{n}/\sigma_{K}(\mathcal{O}_{K}) \rightarrow C_{0}$.}

\Preuve{By lemma~\ref{HcapKtrivial} above, the quotient of $\C^{n}/\sigma_{K}(\mathcal{O}_{K})$ by $H$ is isomorphic to the quotient of $\C^{s+t}=\C^{n}/H$ by $\pi_{H}(\sigma_{K}(\mathcal{O}_{K}))$, which is a Cousin group by lemma~2.4 of \cite{Oeljeklaus:2005aa}. }

\Lemme{\label{intersectionTriviale}The subgroups $H$ and $(\S^{1})^{n}$ of $\C^{n}/\sigma_{K}(\mathcal{O}_{K})$ intersect trivially.}
\Preuve{An element of this intersection corresponds to an element of the linear subspace $H$ satisfying the following equation:
\[
P_{\mathcal{B},\mathcal{B}_{K}}\left(\begin{array}{c}
x_{1}\\
\\
\vdots\\
\phantom{\vdots}\\
\\
x_{n}
\end{array}
\right) = \left(\begin{array}{c}
0\\
\vdots\\
0\\
z_{s+t+1}\\
\vdots\\
z_{s+t+t}
\end{array}
\right),
\]
where $P_{\mathcal{B},\mathcal{B}_{K}}$ is the change-of-basis matrix from the canonical basis to $\mathcal{B}_{K}$ and the $x_{i}$'s are real. 
The $i$-th and $(i+t)$-th components of every vector of the basis $\mathcal{B}_{K}$ are conjugated (for $i=s+1, ..., s+t$). Since the $x_{i}$'s are real numbers, the $i$-th and $(i+t)$-th components of the vector $(0, ..., 0, z_{s+t+1}, ..., z_{s+t+t})$ are also conjugated, hence $z_{s+t+1} = ... = z_{s+t+t} =0$.
}

\Remarque{The quotient of $(\C^{*})^{n}$ by $(\S^{1})^{n}$ is given by the map 
\[\application{\ord}{(\C^{*})^{n}}{\R^{n}}{(z_{1}, ..., z_{n})}{(-\ln|z_{1}|, ..., -\ln|z_{n}|).}\]
With equation~(\ref{injectionCt}), we see that the space $\ord(H) \cong \R^{2t}$ is generated by the $2t$ vectors $\Re(h_{j})$ and $\Im(h_{j})$ (for $j\in\{1,..., t\}$), in other words, $\ord(H)=\widetilde{H}$.}

\subsection{Manifolds with corners}
In \cite{Battisti:2013aa}, one can find the definition of the \emph{manifold with corners} of a (non-necessarily rational) fan $\Delta$ of $\R^{n}$, denoted by $\Mc(\Delta)$. Heuristically, this is a partial compactification of $\R^{n}$ obtained the following way: for every cone $\sigma$ in $\Delta$ one sends a complementary subspace of $\operatorname{span}(\sigma)$ ``at infinity'' in the direction of $\sigma$. One also defines a topology on this space and when $\Delta$ is rational, this space is (homeomorphic to) the manifold with corners of $\Delta$ in the usual toric-geometrical sense, i.e. $X_{\Delta}/(\S^{1})^{n}$.

We need the following two lemmas:

\Lemme{[\cite{Battisti:2013aa}, lemma~1.17]\label{coneCompact}Let $\sigma$ be a cone of a fan $\Delta$. Then the closure $S$ of $\sigma\subset\Mc(\Delta)$ in $\Mc(\Delta)$ is compact.}

\Lemme{[\cite{Battisti:2013aa}, lemma~2.3]\label{equivCond2}Let $\Delta$ be a fan in $\R^{n}$ and $E\cong \R^{k}$ be a linear subspace of $\R^{n}$. The action of $E$ on $\Mc(\Delta)$ is proper if and only if the restriction of the quotient map $\pi: \R^{n} \rightarrow \R^{n}/E$ to the support $|\Delta|$ of $\Delta$ is injective.}

We state the following definition and proposition, the proof of which is straightforward:
\Definition{\label{defActionGpeEv}Let $\Delta$ be a (non-necessarily rational) fan of $\R^{n}$ and $G$ a discrete group of $GL_{n}(\R)$. We say that $G$ \textbf{acts} on $\Delta$ if for every cone $\sigma \in \Delta$ and every $g \in G$, one has $g(\sigma) \in \Delta$. This action is called \textbf{free} if $g(\sigma)\neq\sigma$ when $g$ is not the identity of $G$ and $\sigma\neq\{0\}$ is a cone of $\Delta$. This action is called \textbf{properly discontinuous} if for every $\sigma \in \Delta$, the set $\{g \in G ~|~ (g (\sigma) \cap \sigma)\setminus\{0\} \neq\emptyset \}$ is finite.}

For a fan $\Delta$, we set $|\Delta|^{*} := |\Delta|\setminus\{0\}$ and $|\Delta|^{*c}$ is the complement of this set in $\R^{n}$. One has:

\Proposition{\label{thActionGroupeSurEventail}Let $\Delta$ be a (non-necessarily rational) fan and $G$ a discrete group acting freely and properly discontinuously on $\Delta$. Then this action induces an action of $G$ on $\Mc(\Delta)\setminus |\Delta|^{*c}$ which is also free and properly discontinuous.}

\Remarque{Notice that $\Mc(\Delta)\setminus |\Delta|^{*c}=(\Mc(\Delta)\setminus \R^{n})\cup|\Delta|^{*}$, i.e. this set is the support of $\Delta$ (except $0$) to which we add the components ``at infinity'' of $\Mc(\Delta)$.}

\subsection{A suitable fan}
\noindent In the vector space $E$, we construct a rational fan $\Delta$ with respect to the lattice $\mathcal{O}_{K} \subset E$. First this fan has to be $W$-invariant. Furthermore, the image of its support under the projecting map $\pi_{\widetilde{H}}$ must be included in an open degenerate proper cone $\Omega\subset \pi_{\widetilde{H}}(E)\cong\R^{s}$ invariant under the action of $W$. \emph{Degenerate} means that the closure of the cone $\Omega$ contains a non-trivial linear subspace of $\pi_{\widetilde{H}}(E)$. If these conditions are satisfied, then there is an action of $W$ on the toric manifold associated to this fan.

The projection $\pi_{\widetilde{H}}$ is injective on $\sigma_{K}(\mathcal{O}_{K})$ and $W$-equivariant. Assume there is a fan $\Delta'$ in $\pi_{\widetilde{H}}(E)$ that is generated by elements of $\pi_{\widetilde{H}}(\sigma_{K}(\mathcal{O}_{K}))$ and whose support is such a cone $\Omega$. Take its pullback via the map $\pi_{\widetilde{H}}$. This is the desired fan $\Delta$.

To construct the fan $\Delta'$ in $\pi_{\widetilde{H}}(E)$, we follow the steps of Sankaran in \cite{Sankaran:1989aa} (theorem~2.5). Again, we only give proofs of results when they need modifications, all the others adapt readily.\\

Let $\Omega$ be an open degenerate cone in $\pi_{\widetilde{H}}(E)$ invariant under $W$. We start by describing the maximal vector subspace contained in its closure: 

\Lemme{[\cite{Sankaran:1989aa}, lemma~1.2]\label{ecritureDeOmegaCommeProduit}The vector space $N$ of maximal dimension contained in the closure of $\Omega$ is of the form 
\[N=\{(x_{1}, ..., x_{s}) \in \pi_{\widetilde{H}}(E)~|~x_{i_{1}} = \cdots = x_{i_{s-h}}=0\},\]
with $h>0$.}
\Remarque{Up to renumbering the coordinates of $\pi_{\widetilde{H}}(E)$ and the $\sigma_{i}$'s simultaneously, we may assume that $N$ is the set \[\{(x_{1}, ..., x_{s}) \in \pi_{\widetilde{H}}(E)~|~x_{1} = \cdots = x_{s-h}=0\}.\]}
This allows us to describe $\Omega$ explicitly: 
\Lemme{[\cite{Sankaran:1989aa}, lemma~1.3]\label{descriptionOmega}One has $\Omega = N \times L_{+}$ where $$L = \{x \in\pi_{\widetilde{H}}(E)~|~x_{s-h+1} = ... =  x_{s}=0\}$$ and $$L_{+} =  \{x \in L~|~ \pm x_{i}>0, i\in\{1,..., s-h\} \}.$$}
\noindent Up to composition with a suitable element of $\mathcal{O}_{K}^{*}$, we may assume that $L_{+}$ is the set $\{x\in L~|~x_{i}>0, i\in\{1,...,s-h\}\}$.\\

\noindent Now, we take an integer $h\in\{1,..., s-1\}$ and consider the expression of $\Omega$ given by lemma~\ref{descriptionOmega}. To ensure that $L_{+}/W$ is a compact manifold, we choose the rank $b$ of $W$ to be equal to the dimension $s-h$ of $L$. In particular, we have $1\leqslant b<s$. Remember that given an integer $b \in \{1, ..., s-1\}$, theorem~\ref{thSankExistenceW} guarantees the existence of a subgroup $W$ of rank $b$ of $\mathcal{O}_{K}^{*,+}$ satisfying Assumption~C.

\Proposition{[\cite{Sankaran:1989aa}, proposition~1.4]The action of $W$ on $\Omega$ is free and properly discontinuous.}

\Proposition{[\cite{Sankaran:1989aa}, proposition~2.2]\label{coneEcrase}If $C \subset \Omega \cup \{0\}$ is a nondegenerate closed cone, for all $\eta \in W$ one has: $\bigcap\limits_{a\in\Z} \eta^{a}C \subset L_{+}$.}
\Remarque{\label{preuveLemmeConesEcr}Let's recall the idea of Sankaran's proof: there is a real number $\delta>0$ such that 
\begin{equation}\label{defCDelta}C\subset C_{\delta}:=\{v\in \Omega ~|~\sum_{i=1}^{b}v_{i}^{2}\geqslant \delta \sum_{j>b} v_{j}^{2}\}.
\end{equation}
It is hence enough to establish the result for $C_{\delta}$. This can be done by observing that there exist $N>1$ and $a\in \N$ such that $\eta^{k}C_{\delta} \subset C_{N^{k}\delta}$ for $k>a$. A crucial remark for the next part of the construction is that the constants $a$ and $N$ can be chosen independently of $\eta \in W$. For further details, see the proof of theorem~2.4 in \cite{Sankaran:1989aa}. Assumption~C and the fact that $b<s$ are used in this proof.}

\Theoreme{[\cite{Sankaran:1989aa}, theorem~2.5]\label{thEventailSankaran}There exists an infinite fan $\Delta'$ of $\R^{s}$ stable under the action of $W$, whose support $|\Delta'|$ is $(\Omega\setminus L_{+})\cup\{0\}$, such that all its cones are generated by elements of $\pi_{\widetilde{H}}(\mathcal{O}_{K})$ and $\Delta'/W$ is a finite set of cones.}

\Remarque{In particular, $\Delta'$ is constructed in such a way that the action of $W$ on it is free and properly discontinuous in the sense of definition~\ref{defActionGpeEv}.}

\section{\label{constr}Construction of the new manifolds}
Now we are able to start the construction. It is divided into two steps. First, we choose a group $W$ as in theorem~\ref{thSankExistenceW} and get an infinite fan on which $W$ acts thanks to theorem~\ref{thEventailSankaran}. We take the quotient $X$ of its associated toric manifold by the complex Lie group $H$ and check that we have an action of $W$ on the complex manifold $X$. 

The second step is to find a suitable open subset $U$ of $X$ so that its quotient under the action of $W$ becomes a compact manifold.

\subsection{Step 1}
Let $K$ be a number field, with $s>0$ real and $2t>0$ complex embeddings. Choose an integer ${b\in \{1, ..., s-1\}}$. Theorem~\ref{thSankExistenceW} shows the existence of a subgroup $W$ of $\mathcal{O}_{K}^{*,+}$ of rank $b$ satisfying Assumption~C.\\
Take a fan $\Delta'$ in $\R^{s}$ as in theorem~\ref{thEventailSankaran}. Note $\Delta$ its preimage $\pi_{\widetilde{H}}^{-1}(\Delta')$, which is a rational fan of $E$ with respect to $\mathcal{O}_{K}$ and consider its associated toric manifold denoted by $X_{\Delta}$. The group $H$ acts as a closed subgroup of $(\C^{*})^{n}$ (see equation~(\ref{injectionCt})) on $X_{\Delta}$ and we have the following:

\Lemme{The action of $H$ on $X_{\Delta}$ is free and proper.}
\Preuve{
To check that this action is proper, we first observe that the group $\widetilde{H}=\ord(H) \cong H$ acts properly on $\Mc(\Delta)$. This is a consequence of lemma~\ref{equivCond2} and the fact that the map $\pi_{\widetilde{H}} : E \rightarrow E/\widetilde{H}$ is injective on $|\Delta|$ by construction. Finally, it is clear that the action of $H$ on $X_{\Delta}$ is proper if and only if the action of $\widetilde{H}$ on $\Mc(\Delta)$ is proper because $(\S^{1})^{n}$ is compact. 

The isotropy group of a point $x \in X_{\Delta}$ is a compact subgroup of $H\cong \C^{t}$ and hence trivial.}

The preceding lemma tells us that the quotient $X:=X_{\Delta}/H$ is a complex manifold and since $(\C^{*})^{n}$ is abelian, the following diagram is commutative:

\begin{equation}
\label{diagrammeCommutatifS}
\xymatrix{
X_{\Delta} \ar[rrr]^{\displaystyle (\S^{1})^{n}}\ar[dd]_{\displaystyle  H \cong \C^{t}} & & & \operatorname{\mathcal{M}c}(\Delta) \ar[dd]^{\displaystyle \ord(H)= \widetilde{H}\cong\R^{2t}}\\
& & & \\
X 
 \ar[rrr]^{\displaystyle q}_{\displaystyle p((\S^{1})^{n})\cong (\S^{1})^{n}~~~~~~}& & & \Mc(\pi_{\widetilde{H}}(\Delta))
}\end{equation}

The action of $(\S^{1})^{n}$ on $X_{\Delta}$ descends to $X$ via the group $p((\S^{1})^{n})\cong \S^{1})^{n}$. We note $q: X \rightarrow \operatorname{\mathcal{M}c}(\pi_{\widetilde{H}}(\Delta))$ the quotient map for this action. A detailed proof of the fact that $\Mc(\Delta)/\widetilde{H}$ and $\Mc(\pi_{\widetilde{H}}(\Delta))$ are homeomorphic can be found in section~2.2.2, pp.~558~\&~559 of \cite{Battisti:2013aa}. 


\Lemme{\label{actionZsurQuotient}The action of $W$ on $X_{\Delta}$ descends to an action on $X$ wich we note $(\eta, x)\mapsto\eta\cdot x$ action.}
\Preuve{We shall prove that the action of $W$ normalizes the subgroup $H$ of $(\C^{*})^{n}$.\\

Let $\eta \in W$, $x \in X_{\Delta}$ and $z=(z_{1}, ..., z_{t}) \in H$. First, $\eta$ is equivariant with respect to
\[\application{\tilde{\eta}}{(\C^{*})^{n}}{(\C^{*})^{n}}{(t_{1}, ..., t_{n})}{(t_{1}^{a_{1,1}}\cdots t_{n}^{a_{1,n}}, ..., t_{1}^{a_{i,1}}\cdots t_{n}^{a_{i,n}}, ..., t_{1}^{a_{n,1}}\cdots t_{n}^{a_{n,n}})}\]
where $(a_{i,j})$ is the matrix of the linear map $\eta$ in the basis $\mathcal{B}_{K}$, see for instance \cite{Ewald:1996aa}, theorem~6.4 p.~244. 

\noindent This means that one has $\eta(\iota(z) x) = \tilde{\eta}(\iota(z)) \eta(x)$. The map $\iota$ is defined in equation~(\ref{injectionCt}). Moreover, $\tilde{\eta}(\iota(z)) = \iota(\sigma_{s+t+1}(\eta)z_{1}, ..., \sigma_{s+t+t}(\eta)z_{t})$, hence we have:
\[\eta(\iota(z) x) = \iota(\sigma_{s+t+1}(\eta)z_{1}, ..., \sigma_{s+t+t}(\eta)z_{t}) \eta(x).\]}

\Lemme{The group $W$ acts on $\Mc(\Delta)$ and this action descends to $\Mc(\pi_{\widetilde{H}}(\Delta))$. Moreover, the map $q$ is equivariant with respect to the action of $W$.}
\Preuve{For an element $x\in X_{\Delta}$ (resp. $y \in \Mc(\Delta)$), we note $H.x$, respectively $\widetilde{H}.y$, the $H$-orbit of $x$, respectively the $\widetilde{H}$-orbit of $y$.
The proof of the previous lemma shows that the action of $W$ descends to $\Mc(\Delta)$. The projection map of $E \subset \Mc(\Delta)$ with respect to $\widetilde{H}$ commutes with the action of $W$, see remark~\ref{actionWDescend}, hence there is an action of $W$  on $\Mc(\pi_{\widetilde{H}}(\Delta))$. Note this action again by $(\eta, \widetilde{H}.y)\mapsto \eta \cdot \widetilde{H}.y$.
Now we prove that $q$ is $W$-equivariant.

Let $H.x \in X$ and $\eta \in W$. By lemma~\ref{actionZsurQuotient} and the commutativity of diagram~($\ref{diagrammeCommutatifS})$, we have: $q(\eta \cdot H.x)  = q(H.\eta(x)) = \widetilde{H}.\ord(\eta(x))$. We first prove that $\ord(\eta(x))=\eta(\ord(x))$ on the dense open set $(\C^{*})^{n}$ of $X_{\Delta}$. 

One has:
\begin{eqnarray*}
\ord(x) &=&-\frac{1}{2\pi}\left(\ln|x_{1}|, ..., \ln|x_{n}|\right),
\end{eqnarray*}
and 
\begin{eqnarray*}
\ord(\eta(x)) &=& -\frac{1}{2\pi}(\ln(|x_{1}^{a_{1,1}}\cdots x_{n}^{a_{1,n}}|), ..., \ln(|x_{1}^{a_{n,1}}\cdots x_{n}^{a_{n,n}}|))\\
& = & -\frac{1}{2\pi}\left(\sum_{i=1}^{n} a_{1,i} \ln|x_{i}|, ..., \sum_{i=1}^{n} a_{n,i} \ln|x_{i}|\right)\\
& = & \eta(\ord(x)).
\end{eqnarray*}
By the density of $(\C^{*})^{n}$ in $X_{\Delta}$, this equality holds on $X_{\Delta}$. Hence 
\[q(\eta \cdot H.x) = \widetilde{H}.\ord(\eta(x)) = \widetilde{H}.\eta(\ord(x)) = \eta\cdot(\widetilde{H}.\ord(x)).\]}

\subsection{Step 2}
The last step consists of finding an open subset $U$ of $X$ on which the action of $W$ is free and properly discontinuous. This shows that the quotient $U/W$ is a complex manifold. It turns out that the action of $W$ on the manifold with corners $\Mc(\pi_{\widetilde{H}}(\Delta))$ carries enough informations for our purpose. In this space, we find an open subset on which the action of $W$ is free and properly discontinuous, with a compact fundamental domain. This proves the compactness of the quotient. Then, we pull everything back to $X$ via the map $q$.\\

Let $\eta_{1}, ..., \eta_{b}$ be generators of $W$. Taking their inverses if necessary, we may assume they all satisfy the following \\

\noindent\textbf{Assumption~C$\boldsymbol{^{+}}$}: for every $i\in\{1,...,b\}$ and every nondegenerate closed cone $C\subset \Omega\cup\{0\}$, one has $\bigcap\limits_{a\in\N} \eta_{i}^{a}C \subset L_{+}$. (Remark~\ref{preuveLemmeConesEcr}.) \\

Define the affine subspace
\[H_{0}:=\{(\underbrace{1, ..., 1}_{b \text{ times}}, x_{1}, ..., x_{s-b})~|~x_{1}, ..., x_{s-b}\in \R\} \subset E/\widetilde{H}\]
and set $H_{i}:= \eta_{i}(H_{0})$ for $i=1,...,b$. \label{ensembleBdomFond} Note $B$ the (unbounded) convex envelope of $H_{0}, H_{1}, ..., H_{b}$. 

\Lemme{\label{lemmePavage}Let $B_{b}$ be the image of $B$ by the projection on the $b$ first coordinates of $E/\widetilde{H}$. Then one has $\displaystyle \bigcup_{\eta \in W}\eta(B_{b}) = (\R^{+}_{*})^{b}$ and the union $\displaystyle \bigcup_{\eta \in W\setminus W_{>1}}\eta(B_{b})$ is bounded, where $W_{>1}$ is the set of elements of $W$ having at least one of their $b$ first coordinates greater or equal to $1$.}
\Preuve{
For the first statement, notice that $W$ is a lattice of $\R^{b}$ via the map $w \mapsto (\ln \eta_{1}, ..., \ln \eta_{b})$. The second one comes from the fact that every element $\eta$ of $W \setminus W_{>1}$ satisfies $\eta_{i}<1$ for all $i\in\{1, ..., b\}$. }

\Lemme{The action of $W$ on $~\mathcal{U}\!:=\!\operatorname{\mathcal{M}c}(\pi_{\widetilde{H}}(\Delta))\setminus |\Omega|^{*c}$ is free, properly discontinuous and admits a compact fundamental domain. The same holds for the action of $W$ on the preimage $U\!:=q^{-1}(\mathcal{U})$ in $X$.}
\begin{proof}The fact that the action is free and properly discontinuous is a consequence of proposition~\ref{thActionGroupeSurEventail}.
By theorem~\ref{thEventailSankaran}, there is a finite set of cones $\Sigma:=\{\sigma_{1}, ..., \sigma_{d}\}$ such that $W\cdot \Sigma = \Delta'$. 
We still assume that $\eta_{1}, ..., \eta_{b}$ is a family of generators of $W$ satisfying assumption~C$^{+}$ 
and $B$ is the set defined directly before lemma~\ref{lemmePavage}.
Let us consider the closure $\overline{D}$ in $\operatorname{\mathcal{M}c}(\pi_{\widetilde{H}}(\Delta))$ of the set

\[D:=\underbrace{\left(\bigcup_{\eta \in W^{+}}\eta(|\Sigma|) \cap B\right)}_{D_{1}} \cup \underbrace{\left( |\Sigma| \cap \bigcup_{\eta \in W_{>1}}\eta(B)\right)}_{D_{2}}.\]
Here, $|\Sigma|$ is the union of the cones of $\Sigma$ and $W^{+}$ is the set of elements of $W$ satisfying assumption~$C^{+}$, along with the identity. 
Notice in particular that $W^{+}\subset W_{>1}$. 
Figure~\ref{dessinDomFond} gives a picture of $\overline{D}$ for $s=2$ and $b=1$.~\\
\begin{figure}[!ht]\
\includegraphics[scale=0.75]{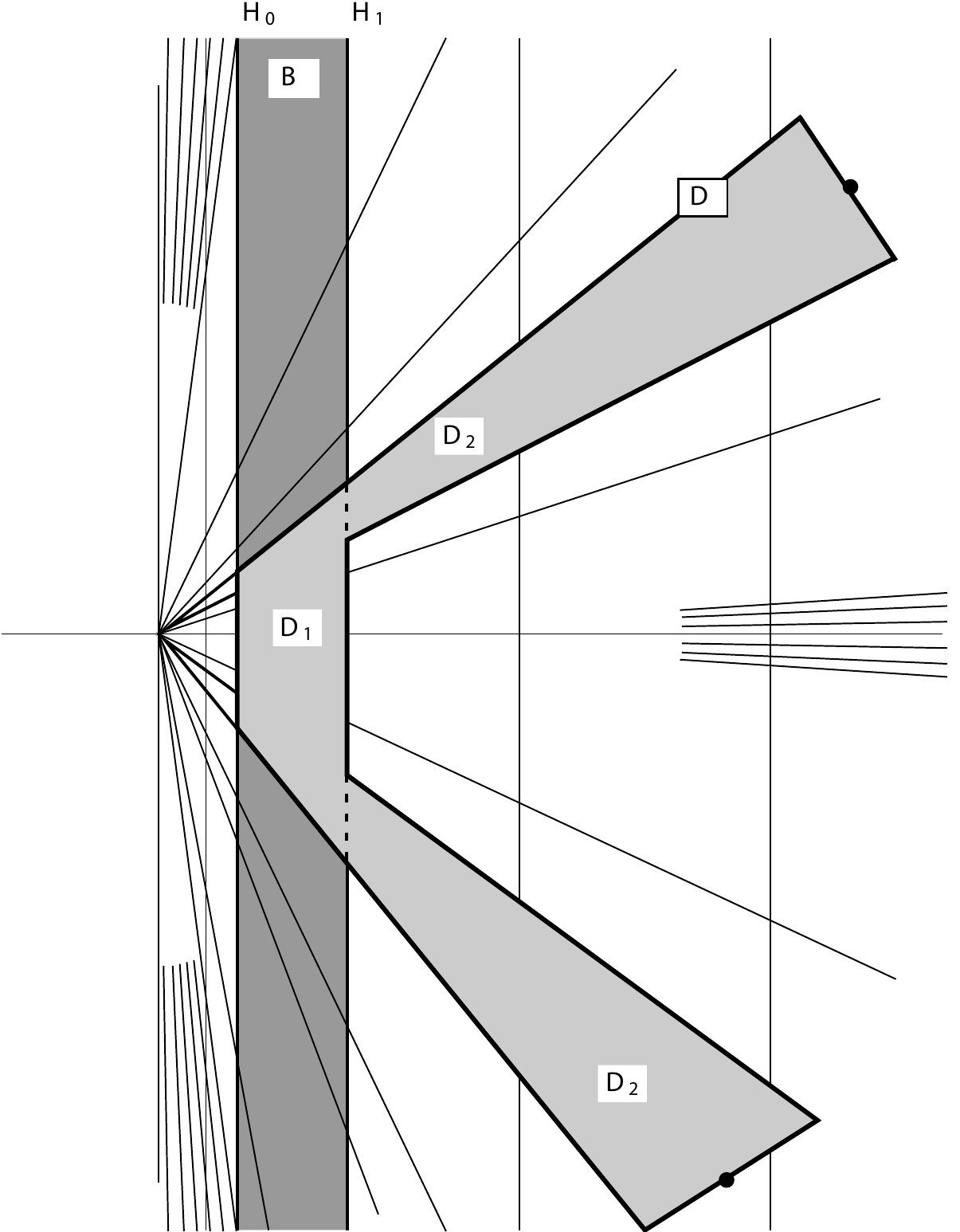}
\caption{The fundamental domain $\overline{D}$.}
\label{dessinDomFond}
\end{figure}

To prove the lemma, it is enough to show that $\overline{D}$ is a compact fundamental domain for the action of $W$.\\
For compactness, first notice that $D_{1}$ is bounded. This is a consequence of remark~\ref{preuveLemmeConesEcr}. Indeed, there is an upper-bound for the $b$ first coordinates of points in $D_{1}$ because of the definition of $B$, while for the $s-b$ last coordinates, we observe that the inclusion $\eta^{a}C_{\delta} \subset C_{N^{a}\delta}$ holds for $a \in \N$ large enough, $\delta > 0$ and $N>1$. Hence the elements of $D_{1}$ lie in a cone $C_{\delta}$ for $\delta>0$ well-chosen. Equation~(\ref{defCDelta}) shows that a subset of $C_{\delta}$ which is bounded in the $b$ first coordinates is also bounded in the $s-b$ last ones.

~\\
Since $\Sigma$ is a finite set, the closure of $D_{2}$ is the union of the closures of $\sigma \cap \bigcup_{\eta \in W_{>1}}\eta(B)$ for $\sigma \in \Sigma$. Such a set is contained in $\sigma \setminus B(0,C)$ for some constant $C> 0$. Indeed, let $x=(x_{1}, ..., x_{s}) \in B$ and $\eta \in W_{>1}$. For all $\eta \in W_{>1}$, there exists an integer $i\leqslant b$ such that $\eta_{i}\geqslant 1$, so without loss of generality we may assume $\eta_{1}\geqslant 1$. The following inequalities hold: 
\[\|\eta(x)\|^{2} = \sum_{i=1}^{s} \eta_{i}^{2} x_{i}^{2} \geqslant \eta_{1}^{2}x_{1}^{2} \geqslant x_{1}^{2} \geqslant \min_{y\in B}y_{1}^{2} > 0.\]
The last inequality is a consequence of the fact that the projection of $B$ on the $b$ first coordinates is a compact set in $\R^{b}$ having no point with a zero coordinate. All the constants $C_{i}:=\min_{y\in B}y_{i}^{2}$ for $i =1, ..., b$ are strictly positive. Let $C$ be the square root of the smallest of these constants. One has $\|\eta(x)\| \geqslant C$. Finally, the closure of $D_{2}$ is compact by lemma~\ref{coneCompact}.\\

In order to finish the proof, we check that $\overline{D}$ is a fundamental domain for the action of $W$. \\
 
First observe that for all $x \in \Omega\setminus L_{+}$, there are two elements $\eta$ and $\eta' \in W$ such that $x \in \eta(|\Sigma|) \cap \eta'(B)$. One has either $\eta^{-1}\eta'\in W^{+} \subset W_{>1}$ or $\eta\eta'^{-1}\in W^{+}$. This respectively implies that either $\eta^{-1}(x) \in |\Sigma| \cap \eta^{-1}\eta'(B)$ or $\eta'^{-1}(x) \in \eta\eta'^{-1}(|\Sigma|)\cap B$. If $x \in L_{+}$, lemma~\ref{lemmePavage} gives the conclusion since $B\cap L_{+}\subset \overline{D}$.

For the remaining case, take a point $x=y+ \infty\cdot \tau$ in $\overline{D}$ on a component at infinity of $\Mc(\pi_{\widetilde{H}}(\Delta))$. The second conclusion of lemma~\ref{lemmePavage} implies then that there is a constant $C'>C$ such that for every cone $\sigma$ of $\pi_{\widetilde{H}}(\Delta)$, one has 
\[\sigma \setminus B(0,C') = \displaystyle \left(\sigma \cap \bigcup_{\eta \in W_{>1}}\eta(B)\right) \setminus B(0,C').\]

Now let $\pi$ be the projection from $\R^{s}$ to $(\R^{s})^{\tau}$ with respect to $L(\tau)$, where $(\R^{s})^{\tau}\oplus L(\tau)=\R^{s}$. Let also $\sigma$ be the $s$-dimensional cone having $\tau$ as a face satisfying $\pi(y) \in \pi(\sigma)$ and let $g$ be an element of the group $W$ such that $g(\sigma) \in \Sigma$. Such a $g$ always exists because of the definition of $\Sigma$. One has: 
\[g(x) = g(y) + \infty \cdot g(\tau).\]
Up to replacing $y$ by $y+w$ with $w\in L(\tau)$, one may assume that $y\in \sigma$, hence $g(y) \in g(\sigma)$ and the conclusion follows. 
We only need to prove the existence of $\sigma$. For this, take a generating set $\{v_{1}, ..., v_{k}\}$ of $\tau$, so that $\tau$ is written as: 
\[\tau = \R_{\geqslant 0}v_{1} + \cdots + \R_{\geqslant 0}v_{k}.\]
Consider the open ball of radius $\varepsilon>0$ centered at $v_{0}:=v_{1}+\cdots + v_{k} \in \tau$. Since the support of $\pi_{\widetilde{H}}(\Delta)$ is $(\Omega\setminus L_{+})\cup\{0\}$, there exists an $\varepsilon>0$ small enough such that $B(v_{0},\varepsilon) \subset \Omega\setminus L_{+}$. Then, by shrinking $\varepsilon$ if necessary, we may assume that $B(v_{0},\varepsilon)$ only intersects the cones containing $\tau$ as a face, in view of the freeness of the $W$-action. Write
\[ B(v_{0},\varepsilon) = B(v_{0},\varepsilon) \cap\!\!\! \displaystyle \bigcup_{\stackrel{\scriptstyle \sigma > \tau}{\dim \sigma = s}} \sigma.\]
The image of $B(v_{0},\varepsilon)$ under the projection $\pi$ is an open neighborhood of $0 \in \pi(\R^{s})$ and the set $\displaystyle \pi(\!\!\bigcup_{\stackrel{\scriptstyle \sigma > \tau}{\dim \sigma = s}}\!\!\sigma)$ is a union of cones, so it is invariant under homothetic transforms of positive ratio. Hence one has $\displaystyle \pi(\!\!\bigcup_{\stackrel{\scriptstyle \sigma > \tau}{\dim \sigma = s}}\!\!\sigma)=\pi(\R^{s})$, there is a cone $\sigma\in\pi_{\widetilde{H}}(\Delta)$ admitting $\tau$ as a face and $\pi(y) \in \pi(\sigma)$.
\end{proof}

\Remarque{In a similar fashion one can describe a compact fundamental domain showing that Sankaran manifolds are compact. Actually, in \cite{Sankaran:1989aa}, p.~47 the author mentions this possibility to obtain compactness, although it did not seem feasible in general.}

Finally, the previous lemma implies:
\Proposition{The quotient $Y\!:=U/W$ is a complex compact manifold of dimension $s+t$.}

\subsection{\label{diviseursSontLVMBGen}Remarks}~\\
\noindent a) The open subset $U\subset X$ intersects the Cousin group $C_{0}$ defined in corollary~\ref{CstGpeCousin}. Note that $U_{0}:=U\cap C_{0}$ is a dense open subset of $U$. The quotient of $U_{0}$ by $W$ is not compact and this is the reason for adding divisors at infinity to $U_{0}$. In fact the irreducible components of these divisors are generalized $\LVMB$ manifolds in the sense of theorem~5.2 of \cite{Battisti:2013aa}.\\
To see this, first observe that the preimage of $U\setminus U_{0}$ in $X_{\Delta}$ is the set $\displaystyle \bigcup_{\sigma\in\Delta\setminus\{0\}} \orb(\sigma)$. Let $\sigma=\R_{\geqslant0}v$ be a $1$-dimensional cone of $\Delta$, with $v\in\sigma_{K}(\mathcal{O}_{K})$. The fan of the toric manifold $\overline{\orb(\sigma)}$, denoted by $\Delta/\sigma$, is obtained by taking the image of every cone of $\Delta$ having $\sigma$ as a face under the projection map with respect to $\operatorname{span}(\sigma)$, see section~VI.4 in \cite{Ewald:1996aa}. Note $\pi_{\sigma}: E\cong\R^{n}\rightarrow E/\operatorname{span}(\sigma)\cong\R^{n-1}$ this projection. One gets a finite fan in $E/\operatorname{span}(\sigma)$ and since $H$ trivially intersects $\sigma_{K}(\mathcal{O}_{K})$, the images of $H$ and $\sigma_{K}(\mathcal{O}_{K})$ under $\pi_{\sigma}$ have a trivial intersection in $E/\operatorname{span}(\sigma)$. Note $\bar\pi_{\widetilde{H}}$ the projection from $E/\operatorname{span}(\sigma)$ with respect to $\pi_{\sigma}(H)$; this is an injective map on $|\Delta/\sigma|$. The fan $\bar\pi_{\widetilde{H}}(\Delta/\sigma)$ is complete because the closure of $\sigma$ in $\Mc(\Delta)$ is compact by lemma~\ref{coneCompact}. All the hypotheses of theorem~5.2 of \cite{Battisti:2013aa} are satisfied.\\

\noindent b) Our construction generalizes other known classes of manifolds. 

i) As for the case $b=0$, which corresponds to the case that $W$ is trivial, one recovers the $\LVMB$ manifolds.

ii) It is also possible to extend the construction to the case $b=s$. In this case the fan $\Delta$ then has to be the trivial fan and this leads to the description of OT manifolds, see~\cite{Oeljeklaus:2005aa} and $Y$ is the compact quotient of $C_{0}$ by $W$.

iii) When replacing the linear subspace $H$ by the trivial space $\{0\}$, the rank $b$ of $W$ varies between $1$ and $s+t-1$ and one gets the construction of Sankaran manifolds as in \cite{Sankaran:1989aa}.

\section{\label{sousSectionInvariants}Invariants and geometric properties}

\Lemme{\label{noHolFunc}The open set $U$ admits no non-constant holomorphic functions.}
\Preuve{The manifold $X$ contains the complex Lie group $C_{0}$ as an open dense subset, see corollary~\ref{CstGpeCousin}. As $U_{0}=U\cap C_{0}$ is an open set of $C_{0}$ stable under the action of the maximal compact subgroup $(\S^{1})^{n}$ of $C_{0}$, it has no non-constant holomorphic functions.}

The following lemma will helps us to compute the Kodaira dimension of our manifolds:
\Lemme{Let $G\cong \C^{t}$ be a closed subgroup of $(\C^{*})^{n}$ and $\Omega$ be an open subset of an $n$-dimensional toric manifold $Z$. Assume that the action of $G$ on $\Omega$ is proper and that the quotient $\Omega/G$ admits no non-constant holomorphic functions. Then the Kodaira dimension of $\Omega/G$ is $-\infty$.}
\Preuve{Let $\pi : \Omega \rightarrow \Omega/G$ be the quotient map. Choose a basis $v_{1}, ..., v_{n}$ of the Lie algebra of $(\C^{*})^{n}$ such that $v_{1}, ..., v_{t}$ is a basis of the Lie algebra of $G$. Define the vector fields $v_{1}^{*}, ..., v_{n}^{*}$ on $\Omega/G$ by $v_{i}^{*}:=d\pi_{y}(v_{i}(y))$ for $y\in \pi^{-1}(x)$ and consider the wedge product $\tau:=v_{t+1}^{*}\wedge ... \wedge v_{n}^{*}\in \Gamma_{\mathcal{O}}(\Omega/G, K^{-1}_{\Omega/G})$. Notice that $\tau$ is non-zero on $\pi((\C^{*})^{n}\cap \Omega)$ while it vanishes on $\pi((Z\setminus(\C^{*})^{n})\cap \Omega)$. If the Kodaira dimension of $\Omega/G$ is not $-\infty$, there is a positive integer $k$ such that there exists a non-trivial section $s\in \Gamma_{\mathcal{O}}(\Omega/G,K^{\otimes k}_{\Omega/G})$. Then $s(\tau^{k})$ is a holomorphic function of $\Omega/G$, hence constant equal to zero, because $\tau^{k}$ only vanishes in the complement of $\pi((\C^{*})^{n}\cap \Omega)$. }

\Proposition{The Kodaira dimension of $Y\!$ is $-\infty$.}
\Preuve{Use the previous lemma with $G=H$, $Z=X_{\Delta}$ and $\Omega=\pi^{-1}(U)$ where $\pi : X_{\Delta} \rightarrow X$ is the quotient map under the action of $H$.}

A consequence of the previous lemma is also:
\Corollaire{The Kodaira dimension of an $\LVMB$ manifold is $-\infty$.}

\Proposition{One has the minoration $\dim H^{1}(Y,\mathcal{O})\geqslant b$.}
\Preuve{Let $\rho: W \rightarrow \C$ be a group homomorphism. We will associate to this homomorphism a principal $\C$-bundle above $Y$ and show that if this bundle is trivial then $\rho$ is also trivial, which will lead to the desired inequality. Consider the action of $W$ on the product $U\times \C$ given by $$\eta.(u,z):= (\eta(u),z+\rho(\eta))$$ for all $\eta \in W$, $u\in U$ and $z\in\C$. The quotient $F:=(U\times \C)/W$ is a principal $\C$-bundle above $Y$. Indeed, the action of $\xi\in\C$ on an element $[u,z]:=W.(u,z) \in F$ is given by $\xi.[u,z]:= [u,z+\xi]$. If $F$ is trivial, it must have a global section, i.e. there is a holomorphic function $f: U \rightarrow \C$ satisfying the equality $f(\eta(u)) = f(u)-\rho(\eta)$ for all $\eta \in W$ and $u\in U$. By lemma~\ref{noHolFunc}, the function $f$ is constant and this implies $\rho \equiv 0$.}

\Corollaire{The fundamental group of the complex manifold $Y$ is isomorphic to $W$, hence one has $b_{1}(Y)=b$. Moreover, $Y$ is non-Kähler.}
\Preuve{As for the fundamental group, the result follows from the fact that the universal covering of $Y$ is simply-connected. By the previous proposition, if $Y$ is Kähler one has $b_{1}(Y) \geqslant 2b$, a contradiction.}

Recall that an element $\alpha$ of $\mathcal{O}_{K}^{*}$ is called a \textbf{reciprocal unit} if $\alpha^{-1}$ is a conjugate of $\alpha$ over $\Q$.
In the following we shall prove that the algebraic dimension of $Y$ is zero if there is such a unit in $W$. 

\Proposition{\label{b2Ouvert}If the group $W$ contains at least one non-reciprocal unit, then $\dim H^{2}(U_{0}/W,\Z)=\binom{b}{2}$.}
\Preuve{The open set $U_{0}/W$ of $Y$ is the quotient of its contractible universal covering $\H^{b}\times \C^{s+t-b}$ by the group $W \ltimes \mathcal{O}_{K}$, hence one has $H^{2}(U/W,\Q)\cong H^{2}(W \ltimes \mathcal{O}_{K},\Q)$. As in \cite{Oeljeklaus:2005aa}, we use the Lyndon-Hochschild-Serre spectral sequence associated to the short exact sequence $$0 \rightarrow \mathcal{O}_{K} \rightarrow W \ltimes \mathcal{O}_{K}\rightarrow W\rightarrow 0.$$ One has $E_{2}^{p,q}=H^{p}(W,H^{q}(\mathcal{O}_{K},\Q))\Rightarrow H^{p+q}(W\ltimes \mathcal{O}_{K},\Q)$ and the following exact sequence:
\begin{eqnarray*}0&\rightarrow& H^{1}(W,\Q^{\mathcal{O}_{K}}) \rightarrow H^{1}(W\ltimes\mathcal{O}_{K},\Q) \rightarrow \overbrace{H^{1}(\mathcal{O}_{K},\Q)^{W}}^{E_{2}^{0,1}} \\
                            & \rightarrow& H^{2}(W,\Q^{\mathcal{O}_{K}})\rightarrow H^{2}(W\ltimes \mathcal{O}_{K},\Q)_{1}\rightarrow \underbrace{H^{1}(W,H^{1}(\mathcal{O}_{K},\Q))}_{E_{2}^{1,1}}
\end{eqnarray*}

where $H^{2}(W\ltimes \mathcal{O}_{K},\Q)_{1}$ is defined by the exact sequence
$$ 0 \rightarrow H^{2}(W\ltimes \mathcal{O}_{K},\Q)_{1} \rightarrow H^{2}(W\ltimes \mathcal{O}_{K},\Q) \rightarrow \underbrace{H^{2}(\mathcal{O}_{K},\Q)^{W}}_{E_{2}^{0,2}}.$$
See for instance \cite{Sah:1974aa}.

If one proves that $E_{2}^{0,1}=E_{2}^{0,2}=E_{2}^{1,1}=0$, the result follows. For the fact that $E_{2}^{0,1}=E_{2}^{1,1}=0$, the proof adapts readily from \cite{Oeljeklaus:2005aa}, proposition~2.3, so we will not repeat it here.\\

As for the group $E_{2}^{0,2}=H^{2}(\mathcal{O}_{K},\Q)^{W} \cong \operatorname{Alt}^{2}(\mathcal{O}_{K},\Q)^{W}$, recall that an element of $\operatorname{Alt}^{2}(\mathcal{O}_{K},\Q)^{W}$ is of the form $\gamma = \sum_{i<j}a_{i,j}\sigma_{i}\wedge \sigma_{j}$ with $a_{i,j} \in \C$. Moreover, the $W$-invariance of $\gamma$ means that for every pair $(i,j)$ such that $a_{i,j}\neq0$, one has $\sigma_{i}(\eta)\sigma_{j}(\eta)=1$ for all $\eta\in W$. Now $W$ contains a non-reciprocal unit $\eta_{0}$ and therefore the relation $\sigma_{i}(\eta_{0})\sigma_{j}(\eta_{0})=1$ can never hold for any choice of $i<j$. Hence $\gamma$ is  trivial and so is the group $E_{2}^{0,2}$.
}


\Proposition{\label{dimAlgNllesVar}If the group $W$ contains at least one non-reciprocal unit, then $U_{0}/W$ contains no complex hypersurface; in particular the algebraic dimension of $Y$ is zero.}
\Preuve{One has the following commutative diagram: 
\begin{equation}
\label{diagrammeOTCousin1}
\xymatrix{
U_{0} \ar[rrr]^{\displaystyle W \cong \Z^{b}} \ar[dd]_{\displaystyle  \R^{s-b}\times(\S^{1})^{n}}^{q} & & & U_{0}/W \ar[dd]^{\displaystyle  \R^{s-b}\times(\S^{1})^{n}}_{q'}\\
& & & \\
(\R_{>0})^{b} \ar[rrr]_{\displaystyle W}& & & (\S^{1})^{b}}
\end{equation}The open set $U_{0}$ is diffeomorphic to $(\R_{>0})^{b}\times\R^{s-b}\times(\S^{1})^{n}$, $W$ acts properly discontinuously on $(\R_{>0})^{b}$ and $U_{0}/W$ is an $\R^{s-b}\times(\S^{1})^{n}$-bundle over $(\S^{1})^{b}$. Because $W$ contains at least one non-reciprocal unit, the previous proposition gives that the map $(q')^{*}: H^{2}((\S^{1})^{b},\Z)\rightarrow H^{2}(U_{0}/W,\Z)$ is injective and the proof of proposition~3.4 of \cite{Battisti:2013ab} adapts readily. Finally, $Y$ only has a finite number of divisors, namely those added at infinity to compactify $U_{0}/W$, hence its algebraic dimension is equal to zero.}

\Remarque{The second Betti number of an OT manifold of simple type was already computed in proposition~2.3 of \cite{Oeljeklaus:2005aa}. In fact the proof of proposition~\ref{b2Ouvert} above shows that we can replace the simplicity condition with the assumption that $s$ is odd.}

\Corollaire{Let $X=X(K,A)$ be an OT manifold. Let $s>0$ and $2t>0$ be the number of real and complex embeddings of $K$ respectively. Then, if $s$ is odd, one has $\dim_{\R} H^{2}(X,\R)=\binom{s}{2}$.}
\Preuve{Notice that this result holds with the same proof as above, under the assumption that the group $A\cong \Z^{s}$ contains at least one non-reciprocal unit. If $s$ is odd, $K$ can not contain any reciprocal unit. Indeed if $K$ contained such a unit, say $\alpha$, the degree of $\Q(\alpha)$ over $\Q$ would be even and divide the degree of $K$.

%
%
%
}

\section{\label{unExempleSalem}An example}
\noindent To conclude this paper, we describe a concrete example of a manifold obtained by the above construction. In what follows we continue to use the same notations. 

\Definition{A \textbf{Salem number} is a real algebraic integer $\gamma>1$ such that all its conjugates are of modulus smaller or equal to $1$, with equality for a least one of them.}

\Remarque{By using the fact that it admits a root of modulus $1$, one easily proves that the minimal polynomial $P$ of a Salem number $\gamma$ is palindromic, i.e. it satisfies $X^{\deg P}P(1/X) = P(X)$. This implies that the roots of $P$ are $\gamma$, $1/\gamma$ and complex numbers of absolute value $1$. In particular, the minimal polynomial of a Salem number has degree at least $4$ and a Salem number is necessarily a unit of even degree.}

\Exemple{\emph{(See \cite{Bertin:1992aa}, p.~85)} It is possible to describe all polynomials of degree $4$ which are minimal polynomials of a Salem number. These are the polynomials with integer coefficients of the form $X^{4}+q_{1}X^{3}+q_{2}X^{2}+q_{1}X+1$ with $2(q_{1}-1)<q_{2}<-2(q_{1}+1)$. The minimal polynomial of the smallest Salem number of degree $4$ is $X^{4} - X^3 - X^{2} -X +1$.}

Now we consider the polynomial $P(X):= X^{4} - X^3 - X^{2} -X +1$. Its roots are $\alpha \approx 1,\!722$ (truncated value), $\alpha^{-1}, \beta$ and $\bar{\beta}$ where $\beta$ is a complex number of modulus $1$ and $\Im (\beta)>0$.
Note $\sigma_{1}, \sigma_{2}, \tau_{1}$ and $\overline{\tau_{1}}$ the associated embeddings of the field $K:=\Q[X]/\langle P \rangle$ in $\R$ (for the first two ones) and $\C$ (for the last two ones).\\ 


One can check by computation that the family $(1, \alpha, \alpha^{2}, \alpha^{3})$ forms an integral basis of $K$, see section~2.6 of \cite{Stewart:2002aa} for a method. The image under the map $\sigma_{K}$ of this family is the following basis of $\C^{4}$: 
\[\mathcal{B}_{K}=\left(
\left(
	\begin{array}{c} 
	1\\
	1 \\ 
	1 \\ 
	1 \\ 
	\end{array}
\right),
\left(
	\begin{array}{c} 
	\alpha\\
	\alpha^{-1} \\ 
	\beta\\ 
	\overline{\beta}\\ 
	\end{array}
\right),
\left(
	\begin{array}{c} 
	\alpha^{2}\\
	\alpha^{-2} \\ 
	\beta^{2}\\ 
	\overline{\beta}^{2}\\ 
	\end{array}
\right),
\left(
	\begin{array}{c} 
	\alpha^{3}\\
	\alpha^{-3} \\ 
	\beta^{3}\\ 
	\overline{\beta}^{3}\\ 
	\end{array}
\right)\right).\]

Moreover, $\alpha$ and $1-\alpha$ are two fundamental units of $\mathcal{O}_{K}^{*}$. Since $s=2$ one has $\Omega=N\times L_{+}$ with $h=\dim N=1 = \dim L=b$ (see lemma~\ref{ecritureDeOmegaCommeProduit} for notations). 

It is clear that the two subgroups $W:=\langle \alpha \rangle_{\Z}$ and $\overline{W}:=\langle 1-\alpha \rangle_{\Z}$ of $\mathcal{O}_{K}^{*,+}$ satisfy Assumption~C and they both are of rank $b=1$. 

The action of $W$ on $\C^{4}$ is given by the following diagonal matrix:
\[ M:=\left(
	\begin{array}{cccc} 
	\alpha & 0 & 0 & 0\\
	0 &\alpha^{-1} &0 &0 \\ 
	0 &0 &\beta &0 \\ 
	0 &0 &0 &\overline{\beta} \\ 
	\end{array}
\right)\!\!.\]

In the basis $\mathcal{B}_{K}$, the matrix of the linear map associated to $M$ is the companion matrix of $P$:
\[C:=\left(
	\begin{array}{cccc} 
	0 & 0 & 0 & -1\\
	1 & 0 & 0 & 1 \\ 
	0 & 1 & 0 & 1 \\ 
	0 & 0 & 1 & 1 \\ 
	\end{array}
\right)\!\!.\]

Here, $H$ is the subgroup $\{(0,0,0,z) \in \C^{4} ~|~ z \in \C\}$ of $\C^{4}$. 
In the basis $\mathcal{B}_{K}$, $H$ is the set $\{(-\beta z,\beta(1-\beta)z,(\overline{\beta} -1)z,z)~|~z\in\C\}$ and the embedding $\iota: H \rightarrow (\C^{*})^{4}$ is given by 
\[\iota(z) = (e^{-2i\pi\beta z},e^{2i\pi\beta(1-\beta)z},e^{2i\pi(\overline{\beta} -1)z},e^{2i\pi z})~\!.\]

The cone $\Omega \subset \R^{2}$ is the open half-plane $\R_{>0}\times \R$. We define a fan of $\R^{2}$ whose cones are generated by elements of $\pi_{\widetilde{H}}(\mathcal{O}_{K})$ and which is invariant under the action of $W$ the following way:  \\

Let $\sigma_{1}$ be the cone positevely generated by the vectors $(1,1)$, $(\alpha,\alpha^{-1})$ and $\sigma_{2}$ the cone generated by $(1,-1)$, $(\alpha,-\alpha^{-1})$ and call $\Delta'$ the fan generated by $W.\{\sigma_{1}, \sigma_{2}\}$. It is clear that $(\Omega\setminus (\R_{>0}\times\{0\})) \cup\{0\}$ is the support of $\Delta'$. For a picture of the situation in $\R^{2}$, we refer to figure~\ref{dessinDomFond}.\\

In this example the divisors that we add to $U_{0}$ so that its quotient by $W$ is compact are chains of Hopf surfaces. 

Indeed, let $\sigma \in \Delta\subset E\cong\R^{4}$ be a $1$-dimensional cone, say $\sigma$ is generated by $(\alpha,\alpha^{-1},\beta,\overline{\beta})$ for simplicity. By the first remark of section~\ref{diviseursSontLVMBGen}, the quotient of $\overline{\orb(\sigma)}$ by $H$ is a generalized $\LVMB$ manifold. The fan $\Delta/\sigma$ of $\R^{3}$ is generated by the images of the two $2$-dimensional cones containing $\sigma$ under $\pi_{\sigma} : \R^{4}\rightarrow\R^{4}/\operatorname{span}(\sigma)$. Up to a linear isomorphism we can assume that this is a subfan of the fan of $\P^{3}(\C)$. In other words, $\overline{\orb(\sigma)}/H$ is an $\LVMB$ manifold. If we adopt the notations of \cite{Lopez-de-Medrano:1997aa}, this situation corresponds to the case $n=4$, $k=3$, $n_{1}=n_{2}=1$ and $n_{3}=2$ hence $\overline{\orb(\sigma)}/H$ is a Hopf surface (section~4.(b), p.~260, ibid.).

Or, one also can notice that $\overline{\orb(\sigma)}/H$ is necessarily a Hopf surface by Potters' theorem on the classification of almost homogeneous complex compact surfaces (\cite{Potters:1969aa}).\\

Hence we have built a complex manifold of dimension $3$ which consists of a dense open set that is the quotient of a Cousin group by a discrete group, compactified by two cycles of Hopf surfaces. A computation shows that each of these Hopf surfaces contains the two elliptic curves $\C/ \langle1, \beta\rangle$ and $\C/ \langle1, 1-\overline{\beta}\rangle$. Now since $1-\overline{\beta} = (\beta-1)/\beta$, these two elliptic curves are isomorphic.\\
 
However, because $W$ only contains reciprocal units, one can not apply proposition~\ref{dimAlgNllesVar} and hence we don't know yet the algebraic dimension of this manifold. On the other hand, if one carries the construction with the group $\overline{W}$ instead of $W$, then we know that the resulting manifold has its algebraic dimension equal to zero.

\bibliographystyle{amsplain}
\bibliography{bibliographie}

\end{document}